\documentclass[11pt]{amsart}
\usepackage{amsfonts,amsmath,amssymb,amsthm,cmtiup} 
\usepackage{mathrsfs}

\usepackage{hyperref}
\usepackage{epigraph}
\begin{document}

\title{Alexander Arhangel'skii} 
\author{Ol'ga Sipacheva}

\begin{abstract}
This is the opening article of the abstract book of conference ``Set-Theoretic Topology 
and Topological Algebra''  in honor of professor 
Alexander Arhangelskii on the occasion of his 80th birthday held in 2018 at Moscow State University.
\end{abstract}

\epigraph{A good mathematician ... never grows old. He remains a child. He remains a dreamer: curious, 
imaginative, free of concrete purpose.}{---\emph{Alexander Arhangel'skii},\\
interview to \emph{The Idler} (1993)}

\maketitle

``Arhangel'skii has a stratospheric reputation and has been regarded over the last thirty 
[written in 2010] years as one of the most important general and set theoretic topologists,''---said 
Peter Collins of Oxford University, United Kingdom, in a letter nominating Arhangel'skii
for the Ohio University Distinguished Professor Award. According to other colleagues from various 
countries, Arhangel'skii is ``one of the foremost general topologists in the world today''; 
his ``prodigious research output is exceptional and proves him to be an original thinker and 
tireless author of top-quality mathematics''; he is a ``great man, creator of $C_p$-theory''; 
he possesses ``almost mystical intuition.''

Alexander Vladimirovich Arhangel'skii was born on March 13, 1938, into a family of artists. His 
father, Vladimir Aleksandrovich Arhangel'skii, was a concert pianist, a student of Rachmaninoff 
and Igumnov and 
close friend of Sofronitsky, an associate professor at the Moscow Conservatory. To be more precise, 
Arhangel'skii Senior had begun as an aeronautical engineer and even become the first elected 
director of Central AeroHydrodynamic Institute, being one of the most promising students of 
Zhukovsky, but his great love for music had made him to change the course of his life. 
Arhangel'skii's mother, Mariya Pavlovna Radimova, was an acknowledged painter, a daughter of the 
celebrated painter and poet Pavel Aleksandrovich Radimov, 
the last chairman of the \emph{Peredvizhniki} (Wanderers) society. 
Strange as it seems at first glance, it may be 
a trait inherited from parents---the ability to tangibly feel harmony and beauty---which is the 
mainstay of his mathematical talent and his incredible intuition. Arhangel'skii once said, 
``Beauty is, for me, a sign of the truth. ... When I think about a mathematical problem or theorem, 
my intuition suggests to me what should be true... And when I try to prove it, I will sometimes 
feel things fitting together harmoniously. It is from this feeling of harmony that I know my way 
will work.''~\cite{karen} In the popular science article~\cite{teaching} on teaching 
mathematics Arhangel'skii concluded several paragraphs of speculation on the harmony and 
beauty of mathematics with the following remarkable words: ``The reader may suspect the author of 
likening mathematics to an art. No, mathematics \emph{is} an art! Imagination, inspiration, and 
illumination are as important for mathematics as for poetry, music, and painting.''

At secondary school Arhangel'skii's favorite subject was literature; he also liked biology and 
physics. Fortunately, in the end, he had decided to study mathematics (although he has always been 
and still remains a bookworm: he reads a lot in Russian, English, French, Italian, Spanish, and, 
very likely, some other languages). So, in 1954, he passed highly competitive entrance examinations and 
entered Moscow State University. During his first year, he took a (required) course in analytical 
geometry taught by Pavel Sergeevich Alexandroff, whose magnetic personality could not but attract 
him, and attended Vitushkin's introductory seminar on function  theory, which provoked his interest 
in set theory and functions of a real variables. These two circumstances have decided his fate: 
when Alexandroff started an introductory seminar in set-theoretic topology, Arhangel'skii attended 
the very first session and decided, right during the session, to specialize in topology. Thus, 
Arhangel'skii the topologist was born in 1955, and Pavel Sergeevich Alexandroff, one of the most 
outstanding Russian mathematicians and one of the foremost pioneers of topology, had become his 
mentor. 

In 1959 Arhangel'skii, still an undergraduate student, obtained his first important result and 
published his first paper, which has become the first of more than five hundred Arhangel'skii's 
papers. That was a metrization theorem for compact spaces based on the new notion of a network in a 
topological space. The result was good, but of most consequence was the introduction of the concept 
of a network, which is now a classical cornerstone concept in general topology. This is very 
typical of Arhangel'skii's work. Although he has proved numerous beautiful theorems, he has always 
been striving to gain insight into the essence of things rather than simply solve existing 
problems. It is hard to find better words than those used by Karen Shenfeld in her 1993 interview 
with Arhangel'skii \cite{karen}: ``He has won for himself an eternal place in the history of 
topology, not so much through the resolution of difficult problems (though he has solved his fair 
share) but rather through the creation of concepts. These concepts ... have become fundamental to 
the ways in which topologists think about space.'' Some of the fundamental notions introduced by 
Arhangel'skii (in addition to the notion of a network) are those of tightness, free sequence, 
strong development, regular base, 
$\tau$-monolithic space, 
$\tau$-balanced, $\tau$-bounded ($\tau$-narrow), and $\tau$-representable topological group,
$\sigma$-paracompactness, $p$-space, 
$\alpha_i$-space, cleavability, Moscow space, and so on; the list might be made very long. 
Arhangel'skii is also very resourceful in inventing amazingly impressive and associative  
names for new objects and properties: feathered, lacy, favorable, Moscow, monolithic, 
radial, cleavable, and so on.

The same year Arhangel'skii graduated, married, and was admitted as a post-graduate student at 
Moscow State University. Two years later, in 1961, Arhangel'skii participated in 
his first international conference, the First Symposium on General Topology in Prague; since then, 
he missed only two TOPOSYMs (because of health problems) and attended many dozens of other 
conferences, mostly as an invited speaker.

In 1962 Arhangel'skii was awarded his candidate of sciences degree (the Russian equivalent of PhD) 
and had to make a hard choice. At that time, his teacher, Alexandroff, headed the Department of 
General Topology at Steklov Mathematical Institute of the Academy of Sciences of the USSR and was 
the chairman of the Department of Higher Geometry and Topology at Moscow State University. 
Arhangel'skii had to decide where to work. He had preferred the University without any hesitation, 
although the Academy was more prestigious and working there did not require teaching. He has always 
though that teaching students and exchanging ideas with them are of great help for research; 
besides, he wanted to be useful. Time has shown that, indeed, Arhangel'skii second vocation is 
teaching. His first PhD student was awarded degree in 1965, and his second and favorite one, 
Mitrofan Choban, in 1969. In all, 37 Arhangel'skii's students have been awarded PhD degree 
(see \cite{genealogy}); many of them  have become doctors of science (a post-doctoral degree, 
Russian analogue of German habilitation) or full professors in different countries. 
Mitrofan Choban was his first student awarded the doctor of science degree. Even now, 
after more than fifty years, when something goes wrong 
with his \emph{Skype}, Arhangel'skii's primary concern is communicating with Mitrofan. 

Arhangel'skii's outstanding teaching ability is closely related to his unique talent for posing 
problems. Virtually all his talks end with dozens of problems, most of which are challenging and 
all interesting. Solving a problem of Arhangel'skii is a good reason for being proud. 
Each term several sessions of his seminar (he halls them ``divertissements'') are devoted entirely 
to problems. During these sessions, all participants are invited to pose and discuss problems, 
but nobody can match Arhangel'skii in this respect. He is very generous in offering his 
problems to students; moreover, he often abstains from thinking on a problem if 
he believes that it may have interesting consequences and can possibly be solved 
by students from the younger generation. He never assigns (but may offer when needed) 
particular problems to his students; instead, he encourages the students to choose 
the problems most interesting to them by themselves. 

Arhangel'skii's problems, which are concerned with very diverse areas, along with his new concepts, 
have greatly influenced the whole development of topology.  Apparently, this talent for posing 
problems originates from Arhangel'skii's  breadth of interest, which has given him a comprehensive 
view of topology and enabled him to see connections hidden from people investigating (maybe, very 
deeply) only certain special things. In his 1991 interview to \emph{Friends of Mathematics 
Newsletter} of Kansas State University \cite{kansas}, Arhangel'skii substantiated his choice of 
topology as follows: ``I consider topology as a fundamental subject. My ideal is to study one of 
the basic concepts, maybe not only of mathematics, but one could also think of it as one of the 
main concepts of a philosophical nature, of life, of common sense,---the concept of continuity. So 
if someone asks me what topology is about, I would say the main idea is to formalize, 
mathematically, the idea of continuity. ... Principles developed in a basic subject like algebra 
and topology help to view mathematics as a whole. They help to unify it. ... With those basic 
principles, you can try to unify mathematics to have a common language. Also, I think, general 
topology is very good for teaching, ... because, being very basic, it does not need other subjects 
for beginners. ... Another thing that is probably related is that there are many, many problems. 
... The problems are at any level: simple, a little bit more difficult, more difficult, and so on. 
And they appear, some new ones. My teacher, Alexandroff, was saying that the most important thing 
in mathematics for young people is to enter the door to creating mathematics; then from that door 
you can move.'' 

In 1966 Arhangel'skii became Doctor of Science. Awarding this highest academic degree in Russia at 
such a young age was (and still remains) quite uncommon, but Arhangel'skii had already been 
internationally recognized for his concept of a network and published a couple dozen papers. One of 
those papers, entitled \emph{Mappings and Spaces} \cite{mappings and spaces}, was a long (about 50 
pages) half-survey, half-research paper, which contained the first systematic reciprocal 
classification of spaces and maps. This paper is also remarkable in that Arhangel'skii showcased 
himself not only as an accomplished mathematician  but also as a poet (look at the epigraph: ``On 
the Edges of Darkness / I sing of Your Galaxies''; this is a loose translation, the Russian 
original reads as something like ``On the branches of Darkness / Blooms the Lilac of Galaxies.'') 

The next year Arhangel'skii applied for membership in the Communist Party of the Soviet Union. 
Joining the Party disagreed with his view that the university should be free from all politics; 
however, in his own words, he simply recognized that, to play a public role in the life of the 
university at that time, he had to become a member. Indeed, many of his students, including the 
author of this note, would never be accepted as post-graduate students without his strong 
``communist'' protection (according to formal rules, only the most active members of Komsomol 
could be accepted, while young mathematicians often considered themselves 
to be too decent or too busy with science for that, or simply  were too careless). However, 
too years had passed before Arhangel'skii was admitted to the Party, because he had signed a letter 
of protest against the psychiatric confinement of his elder colleague Esenin-Volpin, who 
was a notable dissident (that was his third political imprisonment). The letter was broadcasted 
by \emph{Voice of America}, and Esenin-Volpin was released almost immediately thereafter. 
The punishment might have been much more severe than it was (Arhangel'skii was reprimanded 
in oral and written form), but Arhangel'skii did not hesitate to put his signature---he ``wanted 
to state that there were some things that [he] could not accept''~\cite{karen}.

In 1969 Arhangel'skii proved that the cardinality of any first-countable compact Hausdorff space 
does not exceed that of the continuum; thereby, he solved an almost fifty-year-old problem of 
Alexandroff and Urysohn. In fact, Arhangel'skii proved a much more general theorem, namely, that 
$|X| \le 2^{\chi(X)L(X)}$ for any $T_1$-space $X$. His proof essentially used the new 
keystone notion of a free sequence introduced by him in a previous paper; as Richard Hodel 
mentioned in~\cite{hodel}, an important legacy of the theorem was ``the emergence of the closure 
method as a fundamental unifying device in cardinal functions.'' In~\cite{hodel} Hodel formulated 
two criteria for a theorem to be great (the theorem should solve a long-standing problem, and it 
must introduce new techniques and generate new results and open problems) and explains in detail 
why Arhangel'skii's theorem satisfies both requirements. 

The first decade of Arhangel'skii's career as a university faculty member was crowned with winning 
the Lenin Komsomol Prize, a very prestigious state award for young scientists, engineers, and 
artists. It should also not be forgotten that his son was born 
during the same decade (he had already had a baby daughter at the beginning of his career).

The next two decades were as productive as the first. Arhangel'skii spent 1972--75 in 
Pakistan, at the University of Islamabad, as the ``official UNESCO expert on topology'' (still 
retaining the position at the Moscow State University), wrote hundreds of 
seminal papers and participated 
in tens of conferences (in fact, his frequent travels abroad began only with 
\emph{perestroika}). In the early 1990s, the situation in Russia was quite uncertain, and in 1993 
Arhangel'skii accepted professorship at Ohio University. Since then, he spent half time in 
Moscow and half time in Athens, Ohio, every year. In 2003 he was awarded the title of Distinguished 
Professor of Ohio University. However, this did not protect him from being fired, among other 
prominent mathematicians,  from the Moscow State University for  spending too much time abroad. This 
university policy had very sad consequences. The university could not create bearable working 
conditions for faculty members and was deprived of those few worthiest ones who were ready to work 
in unbearable conditions. The university teaching traditions have been lost. Soon after that 
yet another, even worse, misfortune befell Arhangel'skii: he had lost his vision almost completely. 
That was a real calamity for a man who could not live without reading. However, Arhangel'skii has 
never given way to despair. First, he learned to listen to audio books; then, little by little, 
he began to read electronic books and journal articles on a computer (huge white letters on black 
background). Now he again teaches at the Moscow State University (and again combines this job with working 
at another place, this time the Moscow Pedagogical State University), again leads a seminar, 
obtain beautiful results, writes papers (during the past year he published six papers), 
and delivers keynotes at conferences. 

The contribution of Arhangel'skii to topology cannot be overestimated. Not only has he introduced  
fundamental concepts and posed seminal problems; he has also performed a systematic
study of classes of maps in relation to topological properties of spaces; 
created the theory of $p$-spaces; investigated the class of symmetrizable spaces; 
proved metrization theorems; made a dominant contribution to the foundation of $C_p$-theory, the 
theory of free topological groups, and the theory of generalized topological groups (such as 
para-, semi-, and quasi-topological groups); developed the theories 
of relative topological properties, cleavable spaces, and weakly normal spaces. At present, 
he works extensively on topological  homogeneity and remainders of compactifications. 

Arhangel'skii has a highly charismatic personality. He is one of those, very few nowadays, 
old-school professors who has made Russian science. He is a very interesting conversationalist and 
can converse eloquently and competently on any subject---literature, music, politics, poetry, 
history, philosophy... 

It seems that the most appropriate concluding words are those written by Arhangel'skii 
himself~\cite{teaching}:

``Mathematical activity is not for everybody. 

``Working with abstract concepts, dealing with them as the most perfect object of real world (and 
this  attitude is characteristic of any mathematician) are possible only for those who loves the 
subject. ...

Doing mathematics requires loving it.''

\end{document}